\documentclass[12pt]{amsart}
\usepackage{amssymb, hyperref}

\newcommand{\Bgp}{{\Z^\N}}

\long\def\forget#1\forgotten{}
\newcommand{\issuenumber}{30}
\newcommand{\issuemonth}{June}
\newcommand{\issueyear}{2010}

\newcommand{\ed}{
\newpage

\section{Unsolved problems from earlier issues}

\begin{issue}
Is $\binom{\Omega}{\Gamma}=\binom{\Omega}{\Tau}$?
\end{issue}

\begin{issue}
Is $\ufin(\cO,\Omega)=\sfin(\Gamma,\Omega)$?
And if not, does $\ufin(\cO,\Gamma)$ imply
$\sfin(\Gamma,\Omega)$?
\end{issue}

\stepcounter{issue}

\begin{issue}
Does $\sone(\Omega,\Tau)$ imply $\ufin(\Gamma,\Gamma)$?
\end{issue}

\begin{issue}
Is $\fp=\fp^*$? (See the definition of $\fp^*$ in that issue.)
\end{issue}

\begin{issue}
Does there exist (in ZFC) an uncountable set satisfying $\sfin(\cB,\cB)$?
\end{issue}

\stepcounter{issue}

\begin{issue}
Does $X \nin \NON(\cM)$ and $Y\nin\mathsf{D}$ imply that
$X\cup Y\nin \COF(\cM)$?
\end{issue}

\begin{issue}[CH]
Is $\split(\Lambda,\Lambda)$ preserved under finite unions?
\end{issue}

\begin{issue}
Is $\cov(\cM)=\fo$? (See the definition of $\fo$ in that issue.)
\end{issue}

\stepcounter{issue}

\begin{issue}
Could there be a Baire metric space $M$ of weight $\aleph_1$ and a partition
$\mathcal{U}$ of $M$ into $\aleph_1$ meager sets where for each ${\mathcal U}'\subset\mathcal U$,
$\bigcup {\mathcal U}'$ has the Baire property in $M$?
\end{issue}

\stepcounter{issue} 

\begin{issue}
Does there exist (in ZFC) a set of reals $X$ of cardinality $\fd$ such that all
finite powers of $X$ have Menger's property $\sfin(\cO,\cO)$?
\end{issue}

\begin{issue}
Can a Borel non-$\sigma$-compact group be generated by a Hurewicz subspace?
\end{issue}

\begin{issue}[MA]
Is there $X\sbst\bbR$ of cardinality continuum, satisfying $\sone(\BO,\BG)$?
\end{issue}

\begin{issue}[CH]
Is there a totally imperfect $X$ satisfying $\ufin(\cO,\Gamma)$
that can be mapped continuously onto $\Cantor$?
\end{issue}

\begin{issue}[CH]
Is there a Hurewicz $X$ such that $X^2$ is Menger but not Hurewicz?
\end{issue}

\begin{issue}
Does the Pytkeev property of $C_p(X)$ imply that $X$ has Menger's property?
\end{issue}

\begin{issue}
Does every hereditarily Hurewicz space satisfy $\sone(\BG,\BG)$?
\end{issue}

\begin{issue}[CH]
Is there a Rothberger-bounded $G\le\Bgp$ such that $G^2$ is not Menger-bounded?
\end{issue}

\begin{issue}
Let $\cW$ be the van der Waerden ideal.
Are $\cW$-ultrafilters closed under products?
\end{issue}

\begin{issue}
Is the $\delta$-property equivalent to the $\gamma$-property $\binom{\Omega}{\Gamma}$?
\end{issue}

\stepcounter{issue}

\stepcounter{issue}

\general\end{document}}

\newcommand{\Cantor}{{\{0,1\}^\N}}

\newcommand{\fc}{\mathfrak{c}}
\newcommand{\fd}{\mathfrak{d}}
\newcommand{\fp}{\mathfrak{p}}

\newcommand{\NON}{{\mathsf   {NON}}}
\newcommand{\COF}{{\mathsf   {COF}}}

\newcommand{\cM}{\mathcal{M}}

\newcommand{\op}{\operatorname}
\newcommand{\cov}{\mathsf{cov}}

\newcommand{\bbR}{\mathbb{R}}

\newcommand{\EdNote}[1]{\par\medskip\noindent\textbf{#1.}}
\newcommand{\fo}{\mathfrak{od}}

\renewcommand{\split}{\mathsf{Split}}
\newcommand{\bq}{\begin{quote}}
\newcommand{\eq}{\end{quote}}
\newcommand{\cO}{\mathcal{O}}
\newcommand{\cB}{\mathcal{B}}
\newcommand{\BG}{\cB_\Gamma}

\newcommand{\BO}{\cB_\Omega}

\newcommand{\sone}{\mathsf{S}_1}    \newcommand{\sfin}{\mathsf{S}_\mathrm{fin}}

\newcommand{\ufin}{\mathsf{U}_\mathrm{fin}}

\newcommand{\nin}{\not\in}

\newcommand{\cW}{\mathcal{W}}

\newcommand{\N}{\mathbb{N}}

\newcommand{\Z}{\mathbb{Z}}

\newcommand{\sbst}{\subseteq}
\newcommand{\by}[2]{\par\hfill\emph{#1}, #2}

\newcommand{\Tau}{\mathrm{T}}
\newcommand{\CE}{\textsc{CE}}

\newtheorem{thm}{Theorem}[section]
\newcommand{\bthm}{\begin{thm}} \newcommand{\ethm}{\end{thm}}
\newtheorem{prop}[thm]{Proposition}
\newcommand{\bprp}{\begin{prop}} \newcommand{\eprp}{\end{prop}}
\newtheorem{fact}[thm]{Fact}
\newcommand{\bfct}{\begin{fact}} \newcommand{\efct}{\end{fact}}
\newtheorem{prob}[thm]{Problem}
\newcommand{\bprb}{\begin{prob}} \newcommand{\eprb}{\end{prob}}
\newtheorem{lem}[thm]{Lemma}
\newcommand{\blem}{\begin{lem}} \newcommand{\elem}{\end{lem}}
\newtheorem{claim}[thm]{Claim}
\newcommand{\bclm}{\begin{claim}} \newcommand{\eclm}{\end{claim}}
\newtheorem{cor}[thm]{Corollary}
\newcommand{\bcor}{\begin{cor}} \newcommand{\ecor}{\end{cor}}
\newtheorem{conj}[thm]{Conjecture}
\newcommand{\bcnj}{\begin{conj}} \newcommand{\ecnj}{\end{conj}}
\theoremstyle{definition}
\newtheorem{defn}[thm]{Definition}
\newcommand{\bdfn}{\begin{defn}} \newcommand{\edfn}{\end{defn}}
\theoremstyle{remark}
\newtheorem{rem}[thm]{Remark}
\newcommand{\brem}{\begin{rem}} \newcommand{\erem}{\end{rem}}
\newtheorem{cnv}[thm]{Convention}
\newcommand{\bcnv}{\begin{cnv}} \newcommand{\ecnv}{\end{cnv}}
\newtheorem{exam}[thm]{Example}
\newcommand{\bexm}{\begin{exam}} \newcommand{\eexm}{\end{exam}}
\newtheorem{issue}{Issue}

\newcommand{\bpf}{\begin{proof}} \newcommand{\epf}{\end{proof}}
\newcommand{\be}{\begin{enumerate}}
\newcommand{\ee}{\end{enumerate}}
\newcommand{\bi}{\begin{itemize}}
\newcommand{\ei}{\end{itemize}}

\newcommand{\general}{\small\vfill\par\noindent\hrulefill\par
\noindent\textbf{Previous issues.} The previous issues of this
bulletin are available online at\\
\url{http://front.math.ucdavis.edu/search?\&t=\%22SPM+Bulletin\%22}
\\[0.1cm]
\textbf{Contributions.} Announcements, discussions, and open problems should be emailed
to \texttt{tsaban@math.biu.ac.il}\\[0.1cm]
\textbf{Subscription.}
To receive this bulletin (free) to your e-mailbox, e-mail us.
}

\newcommand{\arXivl}[4]{\subsection{#2}{#4}\par\hfill{\arx{#1}}\par\hfill\emph{#3}}
\newcommand{\arXiv}[3]{\subsection{#2}\mbox{}\par\hfill{\arx{#1}}\par\hfill\emph{#3}}

\newcommand{\AMS}[3]{\subsection{#1}\mbox{}\par\hfill{\texttt{#3}}\par\hfill\emph{#2}}
\newcommand{\SPMBul}{\textbf{$\mathcal{SPM}$ Bulletin}}

\newcommand{\arx}[1]{\url{http://arxiv.org/abs/#1}}

\title[$\mathcal{SPM}$ Bulletin \textbf{\issuenumber} (\issuemonth{} \issueyear)]{%
$\mathcal{SPM}$ Bulletin\\[0.5cm]
Issue number \issuenumber: \issuemonth{} \issueyear{} \CE{}}

\begin{document}
\maketitle


\section{Editor's note}

Due to shortage in free time for editing the \SPMBul{},
I am trying for a while a ``slim'' version containing mainly,
when the papers are not in the core of the field of selection principles,
titles and links. By clicking a link, you will get directly
to the webpage of the paper, including the abstract and full paper.
This will make it possible to post fewer issues and
save time for the editor and the readers.

The new style make the table of contents less important, and we therefore
do not include it anymore.

We apologize for not being able to keep up with the earlier format.

\medskip

\by{Boaz Tsaban}{tsaban@math.biu.ac.il}

\hfill \texttt{http://www.cs.biu.ac.il/\~{}tsaban}

\section{Research in SPM}

\noindent\emph{Note.} The division between the present section and the next one
is somewhat artificial, and cannot be precise.

\arXivl{1001.5400}
{Ideals which generalize $(v^0)$}
{Piotr Kalemba and Szymon Plewik}
{We consider ideals $d^0(\mathcal{V})$ which are generalizations of the ideal
$(v^0)$. We formulate couterparts of Hadamard's theorem. Then, adopting the
base tree theorem and applying Kulpa-Szyma\'nski Theorem, we obtain $
cov(d^0(\mathcal{V}))\leq add(d^0(\mathcal{V}))^+$.}

\arXivl{1002.0894}
{Preserving the Lindel\"of property under forcing extensions}
{Masaru Kada}
{We investigate preservation of the Lindel\"of property of topological spaces
under forcing extensions. We give sufficient conditions for a forcing notion to
preserve several strengthenings of the Lindel\"of property, such as
indestructible Lindel\"of property, the Rothberger property and being a
Lindel\"of $P$-space.}

\arXivl{1002.4419}
{Remarks on the preservation of topological covering properties under Cohen forcing}
{Masaru Kada}
{Iwasa investigated the preservation of various covering properties of
topological spaces under Cohen forcing. By improving the argument in Iwasa's
paper, we prove that the Rothberger property, the Menger property and selective
screenability are also preserved under Cohen forcing and forcing with the
measure algebra.}

\arXivl{1002.2883}
{A unified theory of function spaces and hyperspaces: local properties}
{S. Dolecki and F. Mynard}
{Many classically used function space structures (including the topology of
pointwise convergence, the compact-open topology, the Isbell topology and the
continuous convergence) are induced by a hyperspace structure counterpart. This
scheme is used to study local properties of function space structures on
$C(X,\mathbb R)$, such as character, tighntess, fan-tightness, strong
fan-tightness, the Fr{\'e}chet property and some of its variants. Under mild
conditions, local properties of $C(X,\mathbb R)$ at the zero function
correspond to the same property of the associated hyperspace structure at $X$.
The latter is often easy to characterize in terms of covering properties of
$X$. This way, many classical results are recovered or refined, and new results
are obtained. In particular, it is shown that tightness and character coincide
for the continuous convergence on $C(X,\mathbb R)$ and is equal to the
Lindel{\"o}f degree of $X$. As a consequence, if $X$ is consonant, the
tightness of $C(X,\mathbb R)$ for the compact-open topology is equal to the
Lindel{\"o}f degree of $X$.}

\arXivl{1003.0714}
{Generalized Luzin sets}
{Robert Ralowski and Szymon Zeberski}
{In this paper we invastigate the notion of generalized $(I,J)$-Luzin set.
This notion generalize the standard notion of Luzin set and Sierpinski set. We
find set theoretical conditions which imply the existence of generalized $(I,J)$-Luzin set.
We show how to construct large family of pairwise non-equivalent
$(I,J)$-Luzin sets. We find a class of forcings which preserves the property of
being $(I,J)$-Luzin set.}

\arXivl{1003.2308}
{Skeletal maps and I-favorable spaces}
{Andrzej Kucharski and Szymon Plewik}
{It is showed that the class of all compact Hausdorff and $I$-favorable spaces
is adequate for the class of skeletal maps.}

\arXivl{1004.0211}
{Hereditarily Hurewicz spaces and Arhangel'skii sheaf amalgamations}
{Boaz Tsaban, Lyubomyr Zdomskyy}
{A classical theorem of Hurewicz characterizes spaces with the Hurewicz
covering property as those having bounded continuous images in the Baire space.
We give a similar characterization for spaces $X$ which have the Hurewicz
property hereditarily. We proceed to consider the class of Arhangel'skii
$\alpha_1$ spaces, for which every sheaf at a point can be amalgamated in a
natural way. Let $C_p(X)$ denote the space of continuous real-valued functions on
$X$ with the topology of pointwise convergence. Our main result is that $C_p(X)$ is
an $\alpha_1$ space if, and only if, each \emph{Borel} (!) image of $X$ in the Baire space is
bounded. Using this characterization, we solve a variety of problems posed in
the literature concerning spaces of continuous functions.
\par
\medskip
To appear in \emph{Journal of the European Mathematical Society}.}

\EdNote{Note}
{Recently, Lev Bukovsk\'y and Jaroslav \v{S}upina discovered an alternative
proof of the main result of this paper. Their paper is expected to be available before long.}

\arXivl{1004.1879}
{The relation of rapid ultrafilters and $Q$-points to van der Waerden ideal}
{Jana Fla\v{s}kov\'a}
{We point out one of the differences between rapid ultrafilters and $Q$-points:
Rapid ultrafilters may have empty intersection with van der Waerden ideal,
whereas every $Q$-point has a non-empty intersection with van der Waerden ideal.
Assuming Martin's axiom for countable posets we also construct a W-ultrafilter
which is not a $Q$-point.}

\arXivl{1005.0074}
{External characterization of I-favorable spaces}
{Vesko Valov}
{We provide both a spectral and internal characterizations of arbitrary
$\mathrm{I}$-favorable spaces. As a corollary we establish that any perfect
image of an openly generated space is $\mathrm{I}$-favorable. In particular,
any image of compact $\kappa$-metrizable space is $\mathrm{I}$-favorable. The
last statement is a generalization of a result due to P. Daniels, K. Kunen and
H. Zhou. We also generalize a result of Bereznicki\v{i} by proving that every
$\mathrm{I}$-favorable subspace of extremally disconnected space is extremally
disconnected.}

\arXivl{1005.0577}
{The character of topological groups: Shelah's pcf theory and Pontryagin-van Kampen duality}
{Cristina Chis, M. Vincenta Ferrer, Salvador Hernandez, Boaz Tsaban}
{The minimal cardinality of a base at the identity in a topological group $G$,
denoted $\chi(G)$, is one of the major invariants of $G$. A celebrated 1936
result of Birkhoff and (independently) Kakutani asserts that $G$ is metrizable
if, and only if, $\chi(G)$ is countable.

We consider the case where $G$ is the
\emph{dual group} of a metrizable group. Using Pontryagin-van Kampen duality
and pcf theory, we show that also in this case, $\chi(G)$ is well behaved, and
that it is determined by the density and the local density of the base,
metrizable group.

We apply our result to compute the character of free abelian
topological groups, extending a number of results of Nickolas and Tkachenko.

This phenomenon is also reformulated in an inner language, not referring to
duality theory. Here, the compact subsets of quotients by compact subgroups of
$G$ determine its character.

For $G$ dual to a metrizable group, $\chi(G)$ is
especially well behaved in the absence of large cardinals. On the other hand,
when large cardinals are available, some freedom is demonstrated using Cohen's
``forcing'' method, answering a question of Bonanzinga and Matveev.

\EdNote{Please comment!} This is a preliminary version, and none of the authors is
an expert in \emph{all} topics dealt with in the paper. We are very interested
in receiving comments and suggestions from experts in topological groups, duality
theory, and set theory (forcing, large cardinals, pcf theory).}

\arXivl{1005.5542}
{Sequential properties of function spaces with the compact-open topology}
{Gary Gruenhage, Boaz Tsaban, and Lyubomyr Zdomskyy}
{Let $M$ be the countably infinite metric fan.   We show that
$C_k(M,2)$ is sequential and contains a closed copy of Arens space $S_2$.  It follows
that if $X$ is metrizable but not locally compact, then $C_k(X)$ contains a closed copy
of $S_2$, and hence does not have the property AP.

We also show that, for any zero-dimensional Polish space $X$,  $C_k(X,2)$
is sequential if and only if $X$ is either locally compact or the derived set $X'$ is compact.
In the case that $X$ is a non-locally compact Polish space whose derived set is compact,
we show that all spaces $C_k(X, 2)$
are homeomorphic, having the  topology determined by an increasing  sequence of Cantor subspaces, the $n$th one
nowhere dense in the $(n+1)$st.}

\arXivl{1006.1808}
{Pcf theory and cardinal invariants of the reals}
{Lajos Soukup}
{The additivity spectrum $ADD(I)$ of an ideal $I$ is the set of all regular
cardinals $\kappa$ such that there is an increasing chain $\{A_\alpha:\alpha<\kappa\}$
in the ideal $I$ such that the union of the chain is not in $I$.
 We investigate which set $A$ of regular cardinals can be the additivity
spectrum of certain ideals.

 Assume that $I$ is the sigma-ideal generated by the
compact subsets of the Baire space or the ideal of null sets.
 For countable sets we give a full characterization of the additivity spectrum
of $I$: a non-empty countable set $A$ of uncountable regular cardinals can be
$\op{ADD}(I)$ in some c.c.c generic extension iff $A=\op{pcf}(A)$.
}

\section{Related research}\label{RA}

\arXiv{0912.5366}
{Getting more colors}
{Todd Eisworth}

\arXiv{1001.0194}
{A coloring theorem for succesors of singular cardinal}
{Todd Eisworth}

\arXiv{1001.0073}
{Nonmeasurability in Banach spaces}
{Robert Ralowski}

\arXiv{1001.0549}
{Remarks on nonmeasurable unions of big point families}
{Robert Ralowski}

\arXiv{1001.0601}
{Zariski topologies on groups}
{Taras Banakh and Igor Protasov}

\arXiv{1001.0596}
{GO-spaces and Noetherian spectra}
{Authors: David Milovich}

\arXiv{1001.0922}
{Understanding Preservation Theorems, II}
{Chaz Schlindwein}

\arXiv{1001.0908}
{Sequential order under CH}
{Chiara Baldovino}

\arXiv{1001.1175}
{Transfinite Approximation of Hindman's Theorem}
{Mathias Beiglb\"ock and Henry Towsner}

\AMS{Diamond, GCH and weak square}
{Martin Zeman}
{http://www.ams.org/journal-getitem?pii=S0002-9939-10-10192-0}

\arXiv{1001.2819}
{Forcing properties of ideals of closed sets}
{Marcin Sabok, Jindrich Zapletal}

\arXiv{1001.4895}
{Precompact noncompact reflexive abelian groups}
{S. Ardanza-Trevijano, M. J. Chasco, X. Dom\'inguez, M. G. Tkachenko}

\arXiv{1001.4888}
{Multiple gaps}
{Antonio Avil\'es, Stevo Todorcevic}

\arXiv{1002.0347}
{A Combinatorial Proof of the Dense Hindman Theorem}
{Henry Towsner}

\arXiv{1002.1599}
{Common idempotents in compact left topological left semirings}
{Denis I. Saveliev}

\AMS{Diamonds}
{Saharon Shelah}
{http://www.ams.org/journal-getitem?pii=S0002-9939-10-10254-8}

\arXiv{1002.2192}
{More on cardinal invariants of analytic P-ideals}
{Barnab\'as Farkas, Lajos Soukup}

\arXiv{1002.2886}
{When is the {I}sbell topology a group topology?}
{S. Dolecki and F. Mynard}

\arXiv{1002.3089}
{Group topologies coarser than the Isbell topology}
{S. Dolecki, F. Jordan and F. Mynard}

\arXiv{1002.3120}
{Relations that preserve compact filters}
{F. Mynard}

\arXiv{1002.3122}
{Products of compact filters and applications to classical product theorems}
{F. Mynard}

\AMS{New book: Structural Ramsey theory of metric spaces and topological dynamics of isometry groups}
{L. Nguyen Van The}
{http://www.ams.org/journal-getitem?pii=S0065-9266-10-00586-7}

\arXiv{1002.4508}
{A new Lindelof topological group}
{Dusan Repovs, Lyubomyr Zdomskyy}

\arXiv{1002.4456}
{A Model Theoretic Proof of Szemer\'edi's Theorem}
{Henry Towsner}

\arXiv{1002.4730}
{Nondiscrete P-Groups Can be Reflexive}
{Jorge Galindo, Luis Recoder-Nu\~nez and Mikhail Tkachenko}

\arXiv{1003.0918}
{Completely nonmeasurable unions}
{Robert Ralowski and Szymon Zeberski}

\arXiv{1003.2425}
{The Stationary Set Splitting Game}
{Paul Larson and Saharon Shelah}

\arXiv{1003.2479}
{Universally measurable sets in generic extensions}
{Paul Larson, Itay Neeman, and Saharon Shelah}

\arXiv{1003.2496}
{Dense families of countable sets below $\fc$}
{Lajos Soukup}

\arXiv{1003.3189}
{A note on Noetherian type of spaces}
{Lajos Soukup}

\AMS{Uniformizing ladder system colorings and the rectangle refining}
{Teruyuki Yorioka}
{http://www.ams.org/journal-getitem?pii=S0002-9939-10-10330-X}

\arXiv{1003.4024}
{Sections, Selections and Prohorov's Theorem}
{V. Gutev and V. Valov}

\arXiv{1003.4030}
{Minimal functions on the random graph and the product Ramsey theorem}
{Manuel Bodirsky, Michael Pinsker}

\arXiv{1003.5983}
{Complexity of Ramsey null sets}
{Marcin Sabok}

\arXiv{1004.0181}
{Conflict free colorings of (strongly) almost disjoint set-systems}
{Andr\'as Hajnal, Istv\'an Juh\'asz, Lajos Soukup, Zolt\'an Szentmikl\'ossy}

\AMS{Fine asymptotic densities for sets of natural numbers}
{Mauro Di Nasso}
{http://www.ams.org/journal-getitem?pii=S0002-9939-10-10351-7}

\arXiv{1004.2083}
{The character spectrum of $\beta\N$}
{Saharon Shelah}

\arXiv{1004.2172}
{Baire class one colorings and a dichotomy for countable unions of $F_\sigma$ rectangles}
{Dominique Lecomte}

\arXiv{1004.5542}
{Linear ROD subsets of Borel partial orders are countably cofinal in the Solovay model}
{Vladimir Kanovei}

\arXiv{1005.1149}
{The Markov-Zariski topology of an abelian group}
{Dikran Dikranjan, Dmitri Shakhmatov}

\arXiv{1005.2590}
{Some more problems about orderings of ultrafilters}
{Paolo Lipparini}

\arXiv{1005.3528}
{Thin-very tall compact scattered spaces which are hereditarily separable}
{Christina Brech and Piotr Koszmider}

\arXiv{1005.4193}
{The extender algebra and $\Sigma^2_1$-absoluteness}
{Ilijas Farah}

\arXiv{1005.4303 (*cross-listing*)}
{Compactness in Banach space theory - selected problems}
{Antonio Avil\'es and Ond\v{r}ej F.K. Kalenda}

\arXiv{1005.5534}
{Linear ROD subsets of Borel partial orders are countably cofinal in Solovay's model}
{Vladimir Kanovei}

\arXiv{1006.1720}
{Wide scattered spaces and morasses}
{Lajos Soukup}

\ed